\documentclass[11pt]{amsart}

\usepackage{amssymb}
\usepackage{enumitem}
\usepackage{hyperref}
\usepackage{tikz}

\usepackage{fullpage}

\newtheorem{theorem}{Theorem}
\newtheorem{lemma}[theorem]{Lemma}

\newtheorem{corollary}[theorem]{Corollary}

\theoremstyle{definition}
\newtheorem*{ack}{Acknowledgments}

\newtheorem{example}[theorem]{Example}

\DeclareMathOperator{\an}{an}

\DeclareMathOperator{\id}{id}
\let\epsilon\varepsilon

\newcommand*\ratmap{\mathbin{\tikz [baseline=0ex,-latex, dashed, ->] \draw [densely dashed] (0em,0.58ex) -- (1.3em,0.58ex);}}

\title{The Kobayashi pseudometric in the presence of log-terminal singularities}

\author{Finn Bartsch}
\address{Finn Bartsch \\
IMAPP Radboud University Nijmegen \\
PO Box 9010, 6500GL \\
Nijmegen, The Netherlands\\}
\email{f.bartsch@math.ru.nl}



\begin{document}

\begin{abstract}
We show that the Kobayashi pseudometric is well-behaved under resolution of singularities of klt type.
This answers a question of Kamenova and Lehn.
\end{abstract}

\maketitle
\thispagestyle{empty}

For the purposes of this note, a \emph{variety} is an integral separated scheme of finite type over the complex numbers.
We slightly abuse notation and identify a variety $X$ with its associated complex-analytic space $X^{\an}$.
If $X$ is a complex-analytic space, we let $d_X$ denote the Kobayashi pseudometric on $X$; we refer to \cite{KobayashiBook} for the definition and the basic properties of the Kobayashi pseudometric.
If $f \colon X \to Y$ is a morphism of complex-analytic spaces and $d$ is a pseudometric on $Y$, we define a pseudometric $f^*d$ on $X$ via $(f^*d)(x_1,x_2) := d(f(x_1),f(x_2))$.
We say that a variety $X$ has singularities \emph{of klt type} if there is an effective $\mathbb{Q}$-divisor $\Delta$ on $X$ such that the pair $(X, \Delta)$ is klt.
(Of course, any variety with log-terminal singularities satisfies this condition by taking $\Delta = 0$.)
We will prove the following result.

\begin{theorem} \label{mainthm}
Let $X$ be a variety with singularities of klt type and let $p \colon X' \to X$ be a proper birational morphism.
Then $d_{X'} = p^* d_X$.
\end{theorem}

This answers a question of Kamenova and Lehn \cite[Question~3.7.(2)]{KamenovaLehn}.
Note that some assumption on the singularities of $X$ is necessary, as the following examples show.

\begin{example}
Let $C \subseteq \mathbb{P}^2$ be a smooth projective plane curve of degree $\geq 4$ and let $X \subseteq \mathbb{P}^3$ be the cone over $C$.
Then every two points of $X$ can be joined by a chain of rational curves, so that $d_X = 0$.
However, letting $X' \to X$ be the blowup in the vertex of the cone, we see that $X'$ is a $\mathbb{P}^1$-bundle over $C$.
Hence, the Kobayashi pseudometric of $X'$ is the pullback of the Kobayashi pseudometric of $C$ along $X' \to C$ and in particular, $d_{X'} \neq 0$.
\end{example}

\begin{example}[C.~Lehn] \label{example:lehn}
Let $S$ be a fake projective plane, i.e.\ a surface of general type whose Betti numbers are $1,0,1,0,1$.
Then the universal cover of $S$ is a ball in $\mathbb{C}^2$ (see \cite[Section~V.1]{CCS}), so that $d_S$ is a metric on $S$.
Let $S \subseteq \mathbb{P}^n$ be a closed immersion and let $X \subseteq \mathbb{P}^{n+1}$ be the cone over $S$.
As before, we have $d_X = 0$, and if we let $X' \to X$ be the blowup in the vertex of the cone, we again see $d_{X'} \neq 0$.
On the other hand, by Hodge theory, both $H^1(S, \mathcal{O}_S)$ and $H^2(S, \mathcal{O}_S)$ vanish, so that $X$ has rational singularities.
\end{example}

While the second example shows that assuming $X$ to have rational singularities is not enough for Theorem~\ref{mainthm} to hold, the singularities in the above examples are not log-canonical.
It would be interesting to know whether Theorem~\ref{mainthm} can be extended to the log-canonical case.

Our proof of Theorem~\ref{mainthm} relies on the following result due to Demailly--Lempert--Shiffman.
We use the notation $\mathbb{D}$ for the open unit disk in the complex numbers $\mathbb{C}$.

\begin{theorem}[Demailly--Lempert--Shiffman] \label{kobayashimetric_curves}
Let $X$ be a variety.
Then the Ko\-ba\-ya\-shi pseudometric $d_X$ may be defined by only considering those morphisms $\mathbb{D} \to X$ whose image is contained in an algebraic curve $C \subseteq X$.
\end{theorem}
\begin{proof}
If $X$ is quasi-projective, this is \cite[Corollary~1.4]{DLS}.
However, we note that the proof given in \emph{loc.\ cit.}\ is incomplete in the case that $X$ is singular.
A more detailed proof, also covering the non-quasi-projective case, is given in \cite[Section~3]{ThaiDuc}.
\end{proof}

Using Theorem~\ref{kobayashimetric_curves}, we can prove the following simple lemma.
We say that a complex space $X$ has \emph{vanishing Kobayashi pseudometric} if for every two points $x_1, x_2 \in X$, we have $d_X(x_1,x_2) = 0$.

\begin{lemma} \label{simplelemma}
Let $X$ and $Y$ be varieties.
Let $f \colon X \to Y$ be a morphism with connected fibers such that for every $y \in Y$, the fiber $f^{-1}(y)$ has vanishing Kobayashi pseudometric.
Suppose that every curve $C \subseteq Y$ can be lifted along $f$ to a curve $C \subseteq X$.
Then $d_X = f^* d_Y$.
\end{lemma}
\begin{proof}
Let $x_1, x_2 \in X$ be two points.
We have to prove that $d_X(x_1, x_2) = d_Y(f(x_1),f(x_2))$.
The distance-decreasing property of the Kobayashi pseudometric immediately implies that $d_X(x_1, x_2) \geq d_Y(f(x_1),f(x_2))$; so let us show the other inequality.
To do so, let $\epsilon > 0$ be arbitrary and choose a chain of disks, i.e.\ morphisms $(\phi_i \colon \mathbb{D} \to Y)_{i=1}^n$ and points $a_i, b_i \in \mathbb{D}$ such that $\phi_1(a_1) = f(x_1)$, $\phi_n(b_n) = f(x_2)$ and $\phi_i(b_i) = \phi_{i+1}(a_{i+1})$ for all $i = 1, \ldots, n-1$, such that $\sum_{i=1}^n d_{\mathbb{D}}(a_i, b_i) \leq d_Y(f(x_1),f(x_2)) + \epsilon$.
By Theorem~\ref{kobayashimetric_curves}, we may assume that the image of every $\phi_i$ is contained in an algebraic curve $C_i \subseteq Y$.
Since, by assumption, $C_i$ can be lifted along $f$, we can lift the morphisms $\phi_i$ to morphisms $\psi_i \colon \mathbb{D} \to X$.
The points $\psi_i(b_i)$ and $\psi_{i+1}(a_{i+1})$ are then contained in the same fiber of $f$.
Consequently, we have $d_X(\psi_i(b_i),\psi_{i+1}(a_{i+1})) = 0$ for every $i = 1, \ldots, n-1$.
Similarly, we have $d_X(x_1, \psi_1(a_1)) = 0$ and $d_X(x_2, \psi_n(b_n)) = 0$.
Applying the triangle inequality, we hence see that the following inequalities hold.
\begin{align*}
d_X(x_1, x_2)
&\leq d_X(x_1, \psi_1(a_1)) + \sum_{i=1}^n d_X(\psi_i(a_i), \psi_i(b_i)) + \sum_{i=1}^{n-1} d_X(\psi_i(b_i), \psi_{i+1}(a_{i+1})) + d_X(\psi_n(b_n), x_2) \\
&= \sum_{i=1}^n d_X(\psi_i(a_i), \psi_i(b_i)) \leq \sum_{i=1}^n d_{\mathbb{D}}(a_i, b_i) \leq d_Y(f(x_1),f(x_2)) + \epsilon
\end{align*}
Since $\epsilon > 0$ was arbitrary, the claim follows.
\end{proof}

We now give two applications of Lemma~\ref{simplelemma}.
Combining the above lemma with a theorem of Graber--Harris--Starr, we quickly obtain the following result.
Recall that a proper variety $X$ is said to be \emph{rationally connected} if any two points $x, y \in X$ can be joined by a rational curve (i.e., a curve whose normalization is $\mathbb{P}^1$).

\begin{corollary} \label{ratconnected_family}
Let $X$ and $Y$ be varieties.
Let $f \colon X \to Y$ be a proper morphism such that for every $y \in Y$, the fiber $f^{-1}(y)$ is a rationally connected variety.
Then $d_X = f^* d_Y$.
\end{corollary}
\begin{proof}
Since the Kobayashi pseudometric of $\mathbb{P}^1$ vanishes identically, it is clear that a rationally connected variety has vanishing Kobayashi pseudometric.
For every curve $C \subseteq Y$, the base change $X \times_Y C \to C$ is a family of proper, rationally connected varieties.
In particular, its generic fiber $V$ is a proper, rationally connected variety over the function field $\mathbb{C}(C)$ of $C$.
Thus, by applying \cite[Theorem~1.2]{GHS} to (a resolution of) $V$, we see that $V$ admits a $\mathbb{C}(C)$-rational point.
The properness of $X \times_Y C \to C$ implies that this $\mathbb{C}(C)$-rational point of $V$ extends to a section of $X \times_Y C \to C$.
Consequently, every curve $C \subseteq Y$ can be lifted to $X$.
Hence Lemma~\ref{simplelemma} applies and we conclude.
\end{proof}

In \cite[Theorem~9.13]{CampanaFourier}, it is claimed that if $X$ is a smooth projective variety, and $f \colon X \ratmap Y$ is its MRC fibration, then we have $d_X = f^* d_Y$.
We note that our Corollary~\ref{ratconnected_family} does not imply this statement, since not every fiber of the MRC fibration is rationally connected (only a general fiber is).

Using results of Hacon and McKernan on varieties with log-terminal singularities \cite{HaconMcKernan}, we can also easily deduce Theorem~\ref{mainthm}.

\begin{proof}[Proof of Theorem \ref{mainthm}]
By \cite[Corollary~1.5.(1)]{HaconMcKernan}, the fibers of $p$ are rationally chain connected (i.e., any two points can be joined by a chain of rational curves).
In particular, their Kobayashi pseudometric vanishes.
Moreover, by \cite[Corollary~1.7.(2)]{HaconMcKernan} (taking $\Delta$ such that $(X, \Delta)$ is klt, $X=S$, $f = \id_X$ in \emph{loc.\ cit.}), every curve $C \subseteq X$ can be lifted along $p$.
Consequently, we are done by Lemma~\ref{simplelemma}.
\end{proof}

\begin{ack}
I thank Ariyan Javanpeykar for his constant support and many helpful discussions.
I thank Fabio Bernasconi and Christian Lehn for helpful comments, and Christian Lehn for suggesting Example~\ref{example:lehn}.
\end{ack}

\bibliographystyle{alpha}
\bibliography{kobayashi_logterminal}{}
\end{document}